\documentclass[a4paper,12pt]{article}

\usepackage[body={160mm,230mm}]{geometry}
\usepackage{a1}
\usepackage[english]{babel}
\usepackage[latin1]{inputenc} 
\usepackage{graphicx}
\usepackage[dvips]{color}
\usepackage{amsfonts,amssymb}

\def\IS{\mathbb S}

\def\mapright#1#2{\smash{\mathop{\hbox to
0.10cm{\rightarrowfill}}\limits^{#1}_{#2}}}
\def\mapleft#1#2{\smash{\mathop{\hbox to
0.10cm{\leftarrowfill}}\limits^{#1}_{#2}}}

\title{Combinatorial Dehn-Lickorish Twists and Framed Link
Presentations of $3$-Manifolds Revisited}
\author{by Sóstenes Lins}
\date{\today}

\begin{document}

\maketitle

\begin{abstract}
From a pseudo-triangulation with $n$ tetrahedra $T$ of an arbitrary
closed orientable connected $3$-manifold (for short, {\em a
$3D$-space}) $M^3$, we present a gem $J\,'$, inducing $\IS^3$, with
the following characteristics: (a) its number of vertices is $O(n)$;
(b) it has a set of $p$ pairwise disjoint couples of vertices
$\{u_i,v_i\}$, each named {\em a twistor}; (c) in the dual $(J\,
')^\star$ of $J\,'$ a twistor becomes a pair of tetrahedra with an
opposite pair of edges in common, and it is named {\em a hinge}; (d)
in any embedding of $(J\,')^\star \subset \IS^3$, the
$\epsilon$-neighborhood of each hinge is a solid torus; (e) these
$p$ solid tori are pairwise disjoint; (f) each twistor contains the
precise description on how to perform a specific surgery based in a
Denh-Lickorish twist on the solid torus corresponding to it; (g)
performing all these $p$ surgeries (at the level of the dual gems)
we produce a gem $G\,'$ with $|G\,'|=M^3$; (h) in $G\,'$ each such
surgery is accomplished by the interchange of a pair of neighbors in
each pair of vertices: in particular, $|V(G\,')=|V(J\,')|$.

This is a new proof, {\em based on a linear polynomial algorithm},
of the classical Theorem of Wallace (1960) and Lickorish (1962) that
every $3D$-space has a framed link presentation in $\IS^3$ and opens
the way for an algorithmic method to actually obtaining the link by
an $O(n^2)$-algorithm. This is the subject of a companion paper soon
to be released.
\end{abstract}

\section{Motivation}

There exists a rather simple algorithm to go from a framed link
inducing a space to a triangulation of the same space. This was
first done in chapter 11 of~\cite{Kauffman and Lins 1994} via {\em
graph encoded 3-manifolds} or {\em gems}. This algorithm was
improved and incorporated to the computational system {\sf BLINK},
in~\cite{LLins 2007}. Thus to get a blackboard framed link from a
gem is a direct task. However, the contrary, given a gem to find by
a polynomial algorithm a blackboard framed link inducing the same
3D-space is, as far as we know, an untouched problem in the
literature. The reason why it is desirable to have such an algorithm
stems from the fact that the quantum invariants are not computable
from a triangulation based presentation of 3D-spaces. The two
languages, triangulations and blackboard framed links have at
present only a one way translation. This paper starts to fix the
situation by providing a linear algorithm to prove the
Lickorish-Wallace Theorem.

\section{An overview of the Algorithm}

Let $T$ be a pseudo triangulation inducing a $3D$-space $M^3$. From
it we construct the sequence
 $$T \ \mapright{barycentric}{subdivision} \ T\,'\ \mapright{dual}{gem}\
 G_{0}\  \mapright{dipole}{cancelations} \ G_{1}
 \mapright{2-dipole\ creations}{spiked\ cactus}
 \ G_2 \ \mapright{twistors}{spanning\ tree\ X} \ G\
  \ \mapright{localizing}{hinges\ in\ X} \ G\,'. $$

The first passage is just the barycentric subdivision of $T$
producing $T\, '$. The barycentric subdivision of a
$3$-pseudo-complex have its $0$-simplexes naturally colored by the
dimensions it represents, namely, $0,1,2,3$, so that each
tetrahedron receives the four colors. Define the color of a face of
a tetrahedron to be the color of its opposite vertex. The second
passage is given by dualization. The $1$-skeleton of the dual cell
complex $(T')^\star$ of $T'$ is a graph $G_0$ so that its vertices
and edges are $0$-cells and $1$-cells of $(T')^\star$. The edges
inherits the color of its dual triangular $2$-face. This coloring of
edges makes $G_0$ into a {\em gem} and from this colored
$1$-skeleton we can recover $T\, '$, whence $M^3$, see \cite{Lins
1995}. In the remaining passages we get gems $G_1$, $G_2$, $G$ and
$G\,'$, all inducing $M^3$. From gem $G\,'$ we will get the gem
$J\,'$ inducing $\IS^3$ by very simple local moves named
$ji$-twists, which maintains the number of vertices.

There is an operation on gems named $k$-dipole $(k=1,2)$ cancelation
which together with its inverse the $k$-dipole creation is capable
of linking any two gems inducing the same $3$-manifold \cite{Ferri
and Gagliardi 1982}, \cite{Lins and Mulazzani 2006}. The cancelation
of a dipole does not change the induced $3$-manifold and decreases
by two the number of vertices of a gem. Therefore in the passage
$G_0 \longrightarrow G_1$ we simplify the gem so that it has no
$k$-dipoles, $k=1,2$. A gem without $1$-dipoles is called a {\em
crystallization}. In the above sequence $G_1, G_2, G$ are
crystallizations. A $2$-dipole cancelation or creation in a
crystallization yields a crystallization.

The objective to be achieved in $G_2$ and the passage $G_1
\longrightarrow G_2$ become transparents in the language of {\em
thin gems} \cite{Carter and Lins 1999}, here rebaptized {\em knits}.
So, in this passage we start by switching from the gem language to
the {\em knit} language, obtaining a knit $N_1$ from $G_1$. From
$N_1$ we produce a {\em spiked cactus knit} $N_2$ having enough
properties to meet our purposes. Then we switch back to the gem
language by obtaining a gem $G_2$ from $N_2$.  The passage $G_1
\longrightarrow G_2$ is entirely obtained by $2$-dipole creations.
It increases the number of vertices, but as we shall see, $|V(G_2)|
< 6\, |V(G_1)|$. The important aspect of the {\em spiked cactus gem}
$G_2$ (corresponding to a spiked cactus knit) is that a certain
graph $G_2^{jk}$, defined from it, is connected. At this point we
fix a spanning tree $X$ of $G_2^{jk}$, $|X|=p+1$. Each one of the
$p$ edges of $X$ corresponds in $G_2$ to a $ji$-pre-twistor. Again,
all pertinent definitions are given in the next section.

To get the passage $G_2 \longrightarrow G$ we start by replacing
each edge of $X$ which corresponds to a $ji$-pre-twist by an
adequate number of parallel edges, each of which corresponding to a
$j$-twistor. This can be done simply by creating special
$2$-dipoles, the so called {\em double-$8$-moves}. After these moves
all the edges of $X$ correspond in the dual $G^\star$ of $G$ to
pairs of tetrahedra with an opposite pair of vertices in common,
named a {\em hinge}. The interior of these pairs of tetrahedra are
pairwise disjoint. The final passage $G \longrightarrow G\,'$ is the
{\em localization of the hinges} and produce truly pairwise disjoint
solid tori. It is accomplished by a local complexifying move that
replaces a pair of vertices by a fixed configuration of 34 vertices.
See Figure~\ref{fig:localizingHinge}. In this work, as a consequence
of the above passages, we prove the following Theorem:

\begin{theorem} \label{theo:mainTheorem}
There is an algorithm which obtains, given a pseudo-triangulation
$T$ with $n$ tetrahedra, $|T|=M^3$, a gem $J\,'$, having $O(n)$
vertices $|J\,'|=\IS^3$, and a $p$-set of pairs of vertices
$H_\ell=\{u_\ell,v_\ell\}$, $\ell=1,\ldots,p,$ $\ell \neq \ell'
\Rightarrow H_\ell \cap H_{\ell\, '}=\emptyset$ of $J\,'$ such that
the $\epsilon$-neighborhood of $H_\ell^\star$ is a $p$-set of
pairwise disjoint solid tori embedded in $\IS^3$. Moreover, the
$ji$-twists of the $p$ $k$-twistors $H_\ell=\{u_\ell,v_\ell\}$
produce a gem $G\,'$, $|V(G\,')|=|V(J\,')|$ with $|G\,'|=M^3$. The
framing of each component of the link induced by the collection of
solid tori is the linking number of the two boundary component of a
specific cylinder, {\em the strip $s_{uvij}$}, obtained from the
$ji$-twist of the $\ell$-th $j$-twistor.
\end{theorem}

This is a strengthening of the classical result of
Wallace~\cite{Wallace 1960} and Lickorish~\cite{Lickorish 1962} and
its proof relies on a linear algorithm. It opens the way for an
algorithmic method for actually obtaining the link by an
$O(n^2)$-algorithm. This is the subject of a companion
paper~\cite{Lins and Lins 2007} currently under preparation: it
awaits a proper computer implementation.

\section{Twistors, antipoles and their weaker versions}
Let $(i,j,k)$ be a permutation of the non-null of colors of a
bipartite gem $(1,2,3)$. An {\em $i$-twistor} in $G$ is a pair of
vertices of the same class of the bipartition $\{u,v\}$ which are in
the same $0i$- and $jk$-gon, in distinct $0j$-, $ik$-, $0k$- and
$ij$-gons. In the dual $G^\star$ of the gem, an $i$-twistor becomes
a pair of tetrahedra $t_u \cup t_v$ with an opposite pair of edges
in common, namely the pair of edges corresponding to the $0i$- and
$jk$-gons containing both $u$ and $v$. Such an structure is named an
$i$-hinge. If $u$ and $v$ satisfy all the connectivity conditions
but are in distinct class of the bipartition, then $\{u,v\}$ is
called an {\em $i$-antipole}. The {\em $ij$-twist} of an $i$-twistor
is the operation of exchanging the $i$- and $j$-colored neighbors of
$u$ and $v$. A $3$-residue in a gem is a connected component of a
subgraph of $G$ induced by three of the four colors. An $ij$-twist
is an internal operation in the class of gems which does not change
its number of $3$-residues. We also have the following proposition.

\begin{figure}[!h]
   \begin{center}
      \leavevmode
     \includegraphics[width=12cm]{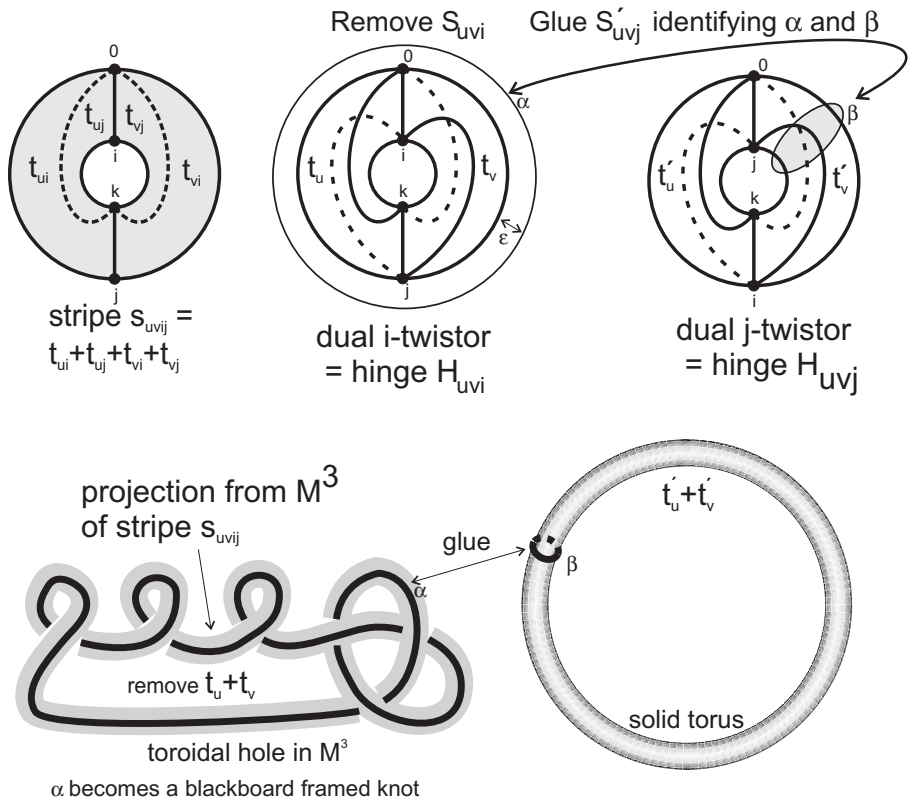}
   \end{center}
   \vspace{-0.7cm}
   \caption{$ij$-surgery on a $i$-hinge: replacement of $t_u
   \cup t_v$, by $t_{{u'}}
   \cup t_{{v'}}$}
   \label{fig:twistorsInTheDual}
\end{figure}

\begin{proposition} \label{twistInTheDual}
For $i,j \in \{1,2,3\}, i\neq j$, in dual terms of hinges an
$ij$-twist of an $i$-twistor in a gem inducing $M^3$ is a
Dehn-Lickorish surgery in $M^3$ \cite{Stillwell 1980}.
\end{proposition}
\begin{proof} We refer to Figure~\ref{fig:twistorsInTheDual}.
For a small $\epsilon$ denote $S_{uvi}$ and $S'_{uvj}$ the solid
tori which are $\epsilon$-neighborhoods of the hinges $t_u \cup t_j$
and $t'_u \cup t'_j$. The two $i$-colored and two $j$-colored faces
of the hinge $t_u \cup t_j$, is topologically a cylinder formed by
four triangles which we call {\em strip} denoting it by
$s_{uvij}=t_{ui} \cup t_{vi} \cup t_{uj} \cup t_{vj}$. Let $\alpha$
be the closed curve in the boundary of $S_{uvi}$ which goes ``just
above'' at an $\epsilon$-distance of the medial curve of $s_{uvij}$.
Let $\beta$ be the boundary of a meridian curve in $S\, '_{uvj}$.
The Dehn-Lickorish surgery is defined by attaching $S\, '_{uvj}$ to
the toroidal hole formed by the removal of $S_{uvi}$ in such a way
as to identify the curves $\alpha$ and $\beta$. Indeed, the only
data needed to perform the surgery is a projection of the curve
$\alpha$ from $M^3$ to a plane with the information of under and
over passes. This projection becomes a blackboard framed knot and
its framing is given by the linking number of the two components of
the the strip $s_{uvij}$. This number can be computed from a general
position projection of the strip in a plane (again keeping the
information of under and over passes).
\end{proof}

For the proof of Theorem~\ref{theo:mainTheorem} we need a weakening
ot the concepts of twistors and antipoles. A {\em $ji$-pre-twistor}
is a pair of vertices in the same class which are in the same $0j$-
and $ik$-gon and in distinct $0i$- and $jk$-gons. If $u$ and $v$
satisfy these connectivity conditions but are in distinct classes,
then the pair $\{u,v\}$ is called a {\em $ji$-pre-antipole}. An {\em
$i$-pair} in a bipartite gem $G$ is a pair of vertices $\{u,v\}$
which are either a $ji$-pre-twistor, a $ki$-pre-twistor, a
$ji$-pre-antipole or a $ki$-pre-antipole. Given a gem $G$ denote by
$G^{jk}$ the graph whose vertices are the $jk$-gons of $G$ and the
edges are the $i$-pairs of $G$. The ends of the edge corresponding
to an $i$-pair $\{u,v\}$ are the vertices corresponding to the
$jk$-gons which contain $u$ and $v$. In particular, $G^{jk}$ may
have parallel edges but not loops.

The only property that we need in the fourth passage, $G_1
\longrightarrow G_2$ is to get by $2$-dipole creations a
crystallization $G_2$ so that $G_2^{jk}$ is connected. We state the
following Conjecture. In it $94$-clusters are the $4$-clusters of
in~\cite{Lins 1995}, where is also defined a rigid gem: a
crystallization whose $3$-residues are $1$-skeletons of polytopes.
Loosely speaking a $94$-cluster consist of a configuration of
$9$-vertices which has a central vertex incident to four square
bigons. A $94$-cluster implies a simplification: the $9$ vertices
become $7$ without changing on the induced $3D$-space. In a rigid
gem each pre-twistor (resp. pre-antipole) is indeed a twistor (resp.
an antipole).
\begin{conjecture}
If $G$ is a rigid $3$-gem without $94$-clusters then $G^{jk}$ is
connected for every choice of distinct $j$ and $k$ in $\{1,2,3\}$.
\end{conjecture}
Unfortunately we have been unable so far to prove this deep
structural property of gems and so we were led to use the {\em
spiked cactus construction}, which follows. This construction is
needed only if $G_1^{jk}$ is not connected, otherwise we take
$G_2=G_1$.

\begin{figure}[!h]
   \begin{center}
      \leavevmode
     \includegraphics[width=15cm]{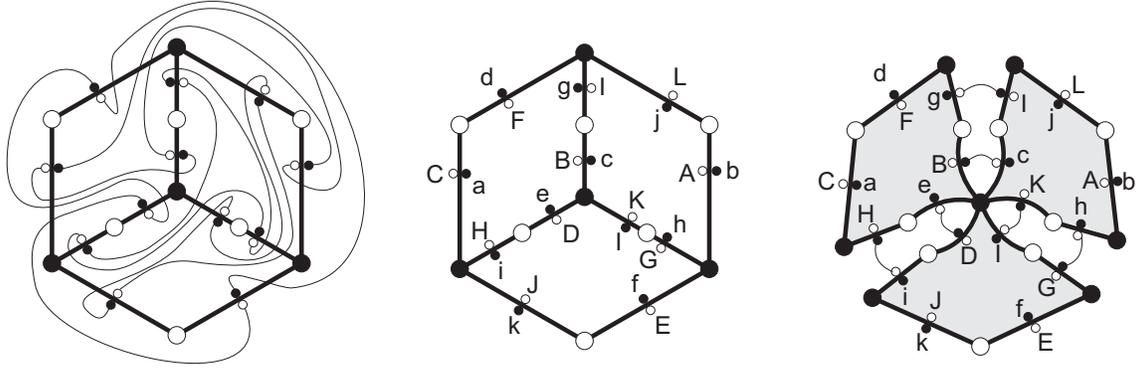}
   \end{center}
   \vspace{-0.7cm}
   \caption{A knit inducing $\IS^1 \times \IS^1 \times \IS^1$ and its ``cactification''}
   \label{fig:knit3torus}
\end{figure}

\section{Knits and spiked cactus knits}

We rename $knit$ to be the object {\em thin gem} introduced for
dimension 3 in~\cite{Lins 1997} and its generalization for dimension
4 given in~\cite{Carter and Lins 1999}. A {\em knit} is a plane
bipartite graph with a perfect match of the sides of the edges. A
$j$-knit $N$ is a knit obtained from a gem $G$ with edge colors
$0,i,j,k$ as follows. Consider the $0$-missing $3$-residues of $G$
embedded in the plane so that the exterior faces is an $ik$-gon.
Choose a representative white vertex to be an interior point of each
$ij$-gon and a representative black vertex to be an interior point
of each $jk$-gon. The black and white vertices are the vertices of a
bipartite graph $N$. Link a black vertex to a white vertex by a
curve so as to cross $G$ only once transversally in the interior of
the $k$-colored edge of $G$ which separates the $ij$-gon and the
$jk$-gon corresponding to the vertices. These curves are the edges
of $N$. Make the cyclic order of these linking edges emanating from
any vertex of $N$ to coincide with the cyclic edge of the dual
$j$-colored edges in the corresponding $ij$-gon or $jk$-gon. Thus
$N$ becomes a graph embedded in the plane.  The vertices of $G$ are
in $1-1$ correspondence with the sides of the edges $N$. These sides
are matched by the $0$-colored edges of $G$ and thus, $N$ is a knit.
A gem $G$ is recoverable from its knit $N$: the angles of the knit
are the $i$- and $k$-colored edges of the gem; the two sides of the
edges of $N$ become the $j$-colored edges of the gem. If $G$ is a
crystallization, $N$ is connected.

A {\em cactus} $j$-knit is one formed by a tree-like arrangement of
polygons and single edges. In a cactus $j$-knit every edge is
incident to the external face (the only white face --- all the
others in light gray). It is always possible to go from any knit
$N_1$ to a cactus knit $N_2'$ by means of a sequence of trivial
angles creations. See the right part of Figure~\ref{fig:knit3torus}.
Each such operation is the creation of a $2$-dipole at the gem
level. We need our knit $N_2$ to have an extra property, namely,
there should not be trivial angles at the black vertices (which
correspond to the $jk$-gons of the associated gem.) This is not the
case of the cactus knit of Figure~\ref{fig:knit3torus}. By extra
trivial angle creations we get easily, in general, a cactus knit
$N_2''$ which meets the extra property.

\begin{figure}[!h]
   \begin{center}
      \leavevmode
     \includegraphics[width=13cm]{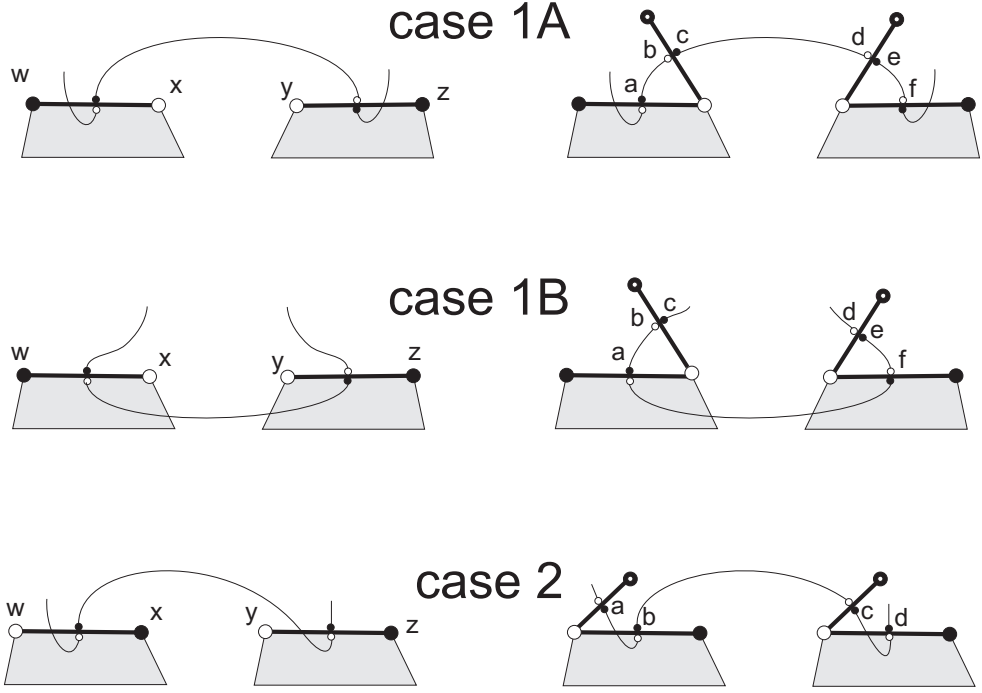}
   \end{center}
   \vspace{-0.7cm}
   \caption{Spikes: $\{a,c\},\{b,d\},\{c,e\},\{d,f\}$  and $\{a,b\},\{b,c\},\{c,d\}$
   are $ji$-pre-twistors }
   \label{fig:spikes}
\end{figure}

To go from $N_2''$ to our desired spiked cactus knit $N_2$ we create
two spikes (two pendant edges) ``trisecting'' each $0$-edge. There
are three cases on how this must be done and they are depicted on
Figure~\ref{fig:spikes}. The creation of each spike is a $2$-dipole
creation in the gem language. In Figure~\ref{fig:spikes} the black
vertices with a white spot correspond to black monovalent vertices
corresponding, in the gem, to $jk$-gons with two vertices. In cases
1A and 1B we might without loss of generality suppose that by going
around the boundary of the external face and in the clockwise
direction the directed edge $wx$ is at the left of the directed edge
$yz$ and that the vertices $x$ and $y$ are white (they might
coincide). Otherwise, they are black and $yz$ is at the left of $wx$
and $a$ and $w$ are white (in the complementary circular path).
Interchanging the labels of $(w,y)$ and $(x,z)$ we get the
assumption holding. For case 2 there is no changing in the argument
if left and right are interchanged. Here is the crucial property of
spiked cactus knits.

\begin{lemma}
Let $N_2$ be a spiked cactus knit of a crystallization $G_2$. Then
every pair $(u,v)$ of edge sides in $N_2$ corresponding to a
$0$-colored edge $(u,v)$ in $G_2$ implies an edge of $G_2^{jk}$
linking the vertices corresponding to the $jk$-gons containing the
vertices $u$ and $v$ of $G$. As a consequence, since $G_2$ is a
crystallization, $G_2^{jk}$ is connected.
\end{lemma} \label{lem:connection}
\begin{proof} The proof follows from the facts that $\{a,c\},\{b,d\},\{c,e\},\{d,f\}$
are $ji$-pre-twistors in the first two cases of
Figure~\ref{fig:spikes} and the same is true for
$\{a,b\},\{b,c\},\{c,d\}$ in the third case. The connectivity
conditions making these pairs pre-twistors are easily checked in the
knit language.
\end{proof}

In consequence of Lemma~\ref{lem:connection}, $G_2$ induces a
connected $G_2^{jk}$. Suppose that $|V(G_2^{jk})|$=p+1 and let $X$
be a subset of $p$ edges of $V(G_2^{jk})$ forming a spanning tree.

\section{Pairwise disjoint solid tori corresponding to $X$}
\begin{figure}[!h]
   \begin{center}
      \leavevmode
     \includegraphics[width=8cm]{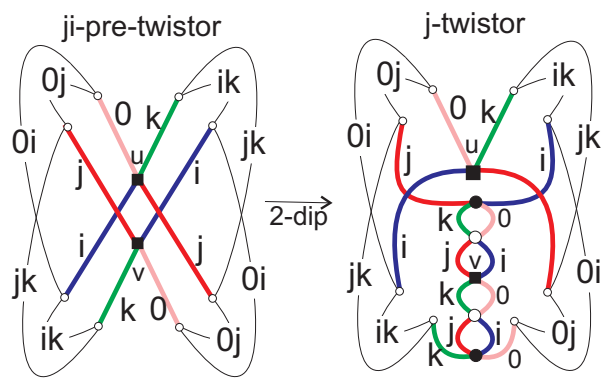}
   \end{center}
   \vspace{-0.7cm}
   \caption{Getting a $j$-twistor from a $ji$-pre-twistor by a $k$-double-8-move at $v$ }
   \label{fig:twistorFromPreTwistor}
\end{figure}

We turn back to the language of gems. The intial step in the passage
$G_2 \longrightarrow G$ is very simple. For each $ji$-pre-twistor in
$X$ which is not a $j$-twistor we fix the situation by creating two
$2$-dipoles near $v$ (a {\em double-8-move}), as depicted in
Figure~\ref{fig:twistorFromPreTwistor}. Note that after this move
the vertices $u$ and $v$ are in distinct $0k$-gons and $ij$-gons.
Thus $\{u,v\}$ becomes a $j$-twistor. A double-8-move replaces a
twistor edge in $G$ by three twistor edges in parallel. Recall that
each edge in $G^{jk}$ is labeled by a pair of vertices forming a
pre-twistor or a pre-antipole. Actually the edges in $X$ correspond
to pre-twistors. By repeating enough double-8-moves at appropriate
vertices we might suppose that in gem $G$ the $p$ edges of $X$ have
distinct $2p$ labels and that each pair of labels is a $j$-twistor.

\begin{figure}[!h]
   \begin{center}
      \leavevmode
     \includegraphics[width=7cm]{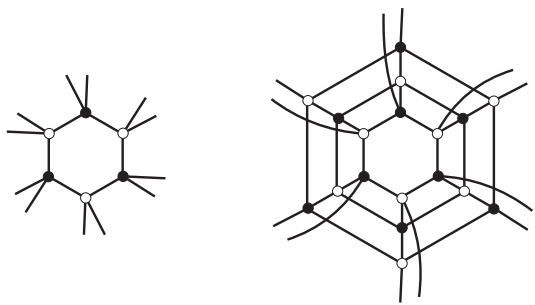}
   \end{center}
   \vspace{-0.7cm}
   \caption{Trisecting a bi-colored polygon in a gem}
   \label{fig:bigonTrisection}
\end{figure}

\begin{figure}[!h]
   \begin{center}
      \leavevmode
     \includegraphics[width=14cm]{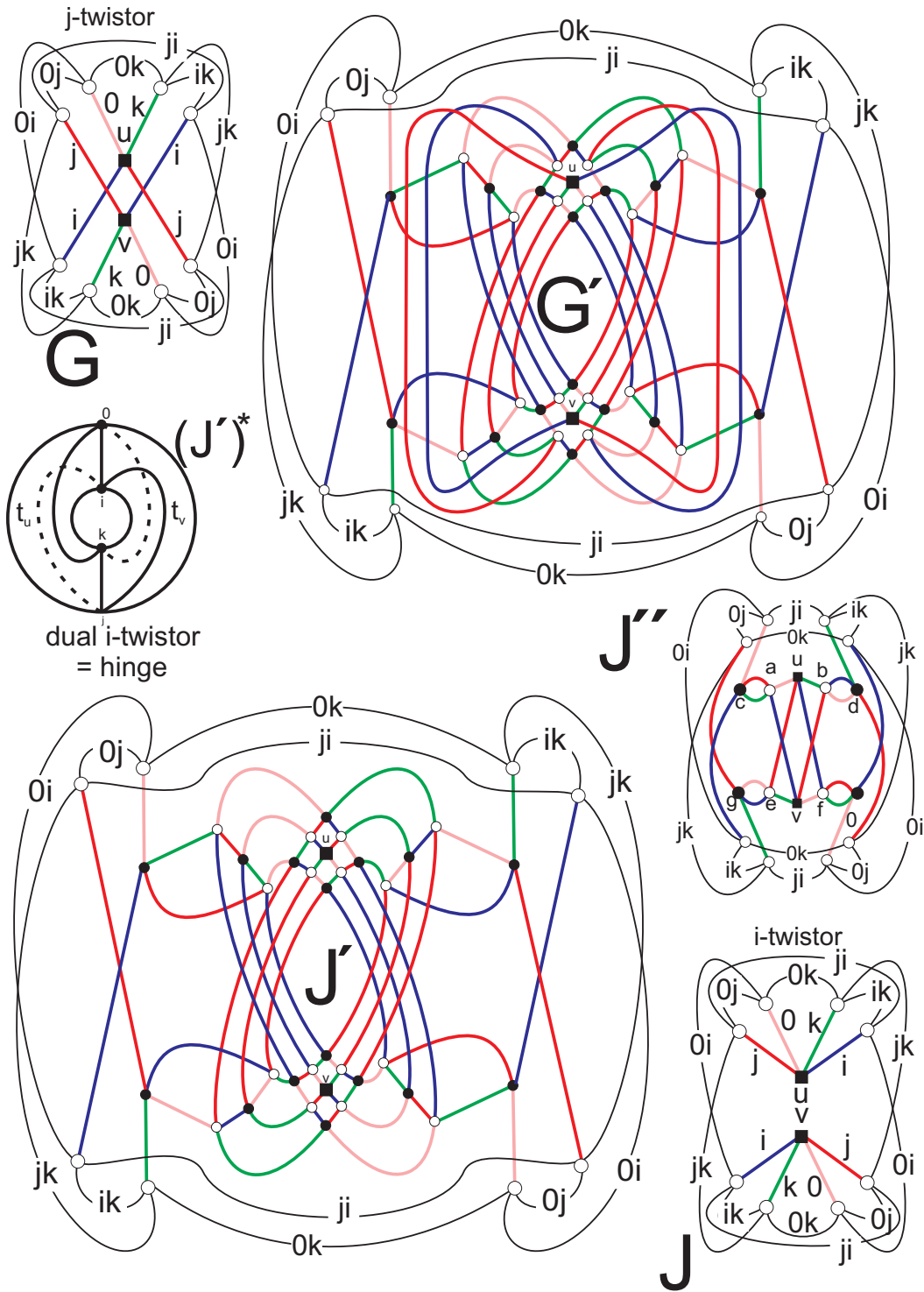}
   \end{center}
   \vspace{-0.7cm}
   \caption{$G\, ' \longrightarrow G$: localizing the hinge}
   \label{fig:localizingHinge}
\end{figure}

From now on we refer to Figure~\ref{fig:localizingHinge}. In the
passage $G \longrightarrow G\,'$ we effect the process of {\em
localizing the hinges}. This means that in $(G\,'    )^\star$ the
set of $p$ pairs of tetrahedra corresponding to $X$ are pairwise
(entirely) disjoint. The passage $G \longrightarrow G\,'$ is
effected ``at the spanning tree $X$ of twistors'' replacing each
pair of vertices of each twistor by the configuration of 34 vertices
depicted in the upper right part of the figure. The localization
moves transform each twistor in $X$ into a local object in the sense
that any order to perform the twists in $X$ produce only gems and
arrive at gem $J$: the $ji$-twistor $\{u_\ell,v_\ell\}$ does not
disturb the other $j$-twistors $\{u_{\ell\, '},v_{\ell\, '}\}$,
$\ell\, ' \neq \ell$. This ``$2\ by\ 34$ replacement'' is essential
for our work. The detailed proof that $G\,'$ simplifies to $G$ by
dipole cancelations is given in the Appendix.

Gem $G\,'$ induces $M^3$ and we prove that performing the
$ji$-twists at the whole set $X$ we get a gem $J\,'$ which induces
$\IS^3$. We observe that to go from $J\,'$ to $G\,'$ we perform the
inverse operation, namely, the $ij$-twist at $X$. It is rather easy
to prove that gems $J\,'\,'$ and $J$ (see
Figure~\ref{fig:localizingHinge}) induce the same space: from $J$ to
$J\,'\,'$ we have four $2$-dipole creations. From $J\,'\,'$ to
$J\,'$ we trisect (see Figure~\ref{fig:bigonTrisection}) the
$0i$-gon and the $jk$-gon incident to $u$ and $v$. These moves are
factorable as $1$-dipoles and $2$-dipole creations. Finally, it is
straightforward to prove that $J$ induces $\IS^3$ because it is a
crystallization having a unique $jk$-gon: from $G$ we can arrive to
$J$ directly by $p$ $ji$-twists and each one of these decreases by
one the number of $jk$-gons. In any crystallization the number of
$0i$-gons and of $jk$-gons coincide.

\section{Appendix: proof that $G\,' \longrightarrow G$ by dipole moves}

\begin{figure}[!h]
   \begin{center}
      \leavevmode
     \includegraphics[width=16cm]{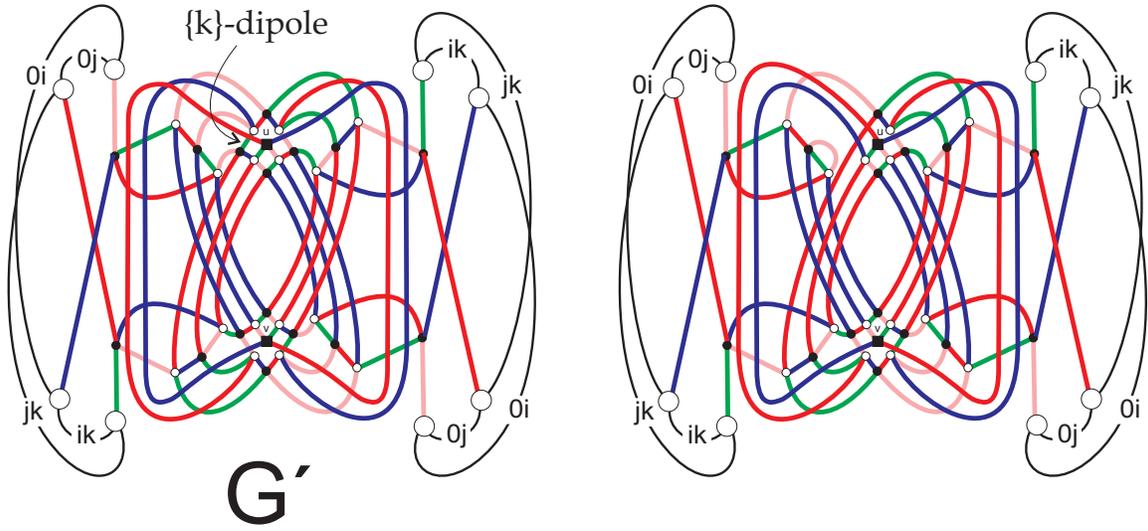}
   \end{center}
   \vspace{-0.7cm}
   \caption{Getting $G$ from $G\,'$: (1/12) canceling a $\{k\}$-dipole}
   \label{fig:CD1}
\end{figure}

\begin{figure}[!h]
   \begin{center}
      \leavevmode
     \includegraphics[width=16cm]{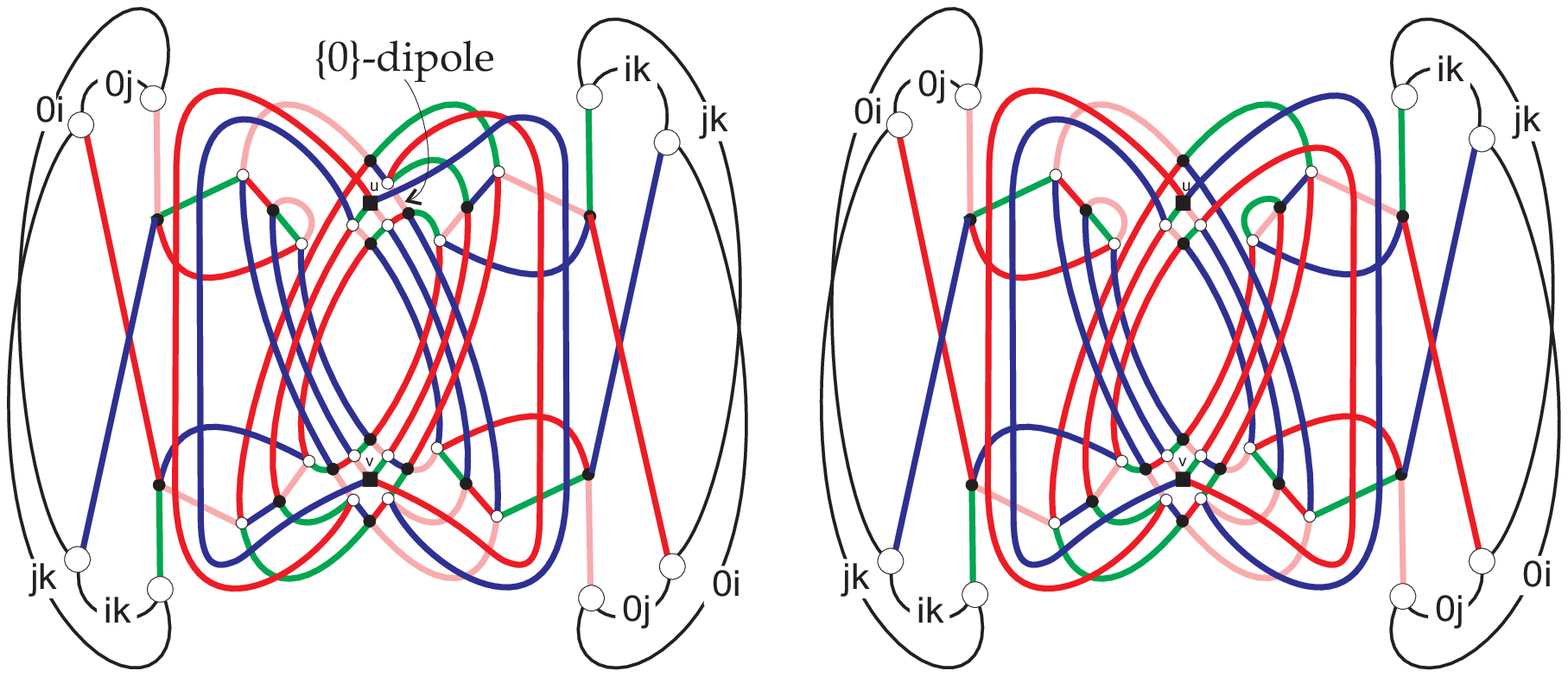}
   \end{center}
   \vspace{-0.7cm}
   \caption{Getting $G$ from $G\,'$: (2/12) canceling a $\{0\}$-dipole}
   \label{fig:CD2}
\end{figure}

\begin{figure}[!h]
   \begin{center}
      \leavevmode
     \includegraphics[width=16cm]{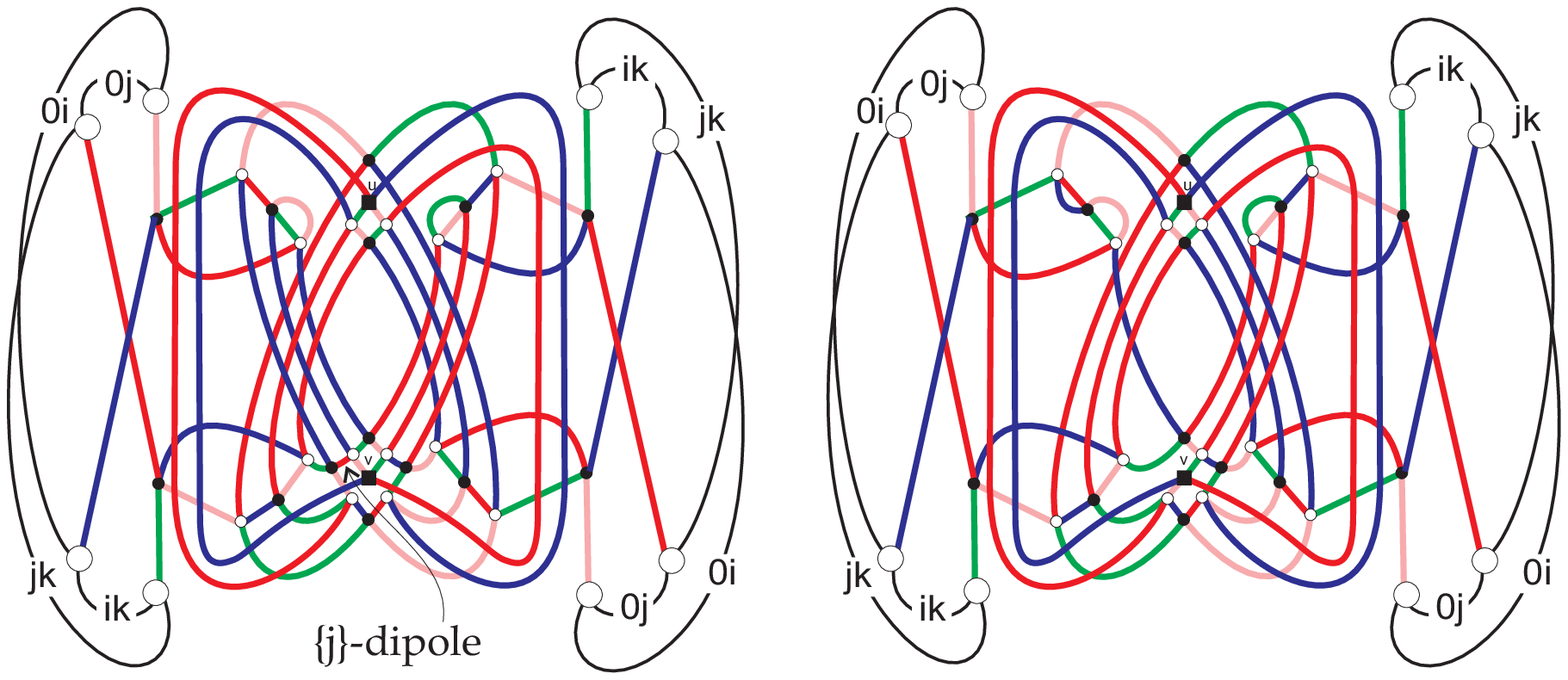}
   \end{center}
   \vspace{-0.7cm}
   \caption{Getting $G$ from $G\,'$: (3/12) canceling a $\{j\}$-dipole}
   \label{fig:CD3}
\end{figure}

\begin{figure}[!h]
   \begin{center}
      \leavevmode
      \includegraphics[width=16cm]{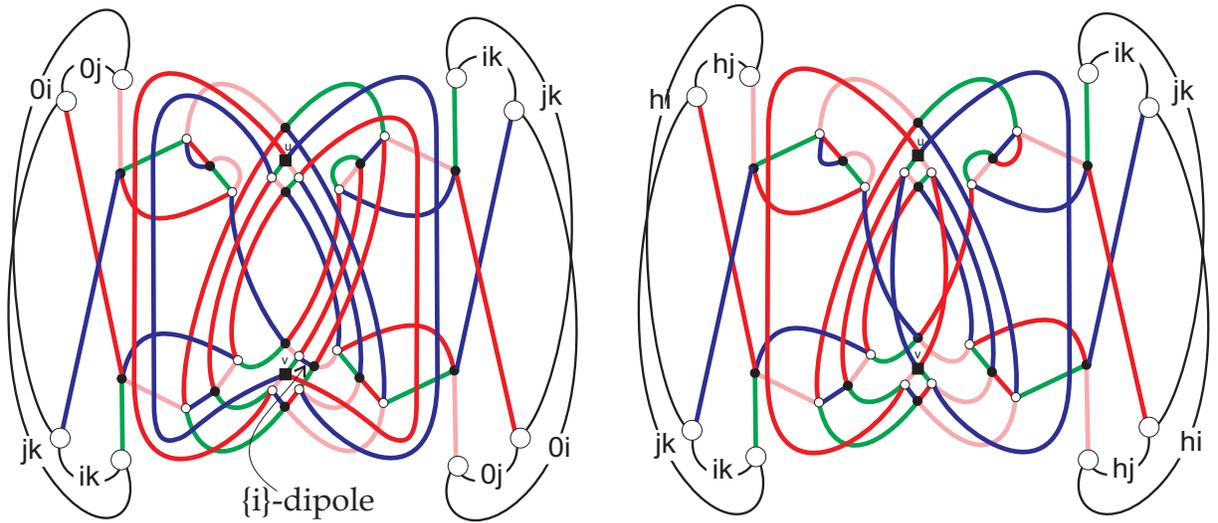}
   \end{center}
   \vspace{-0.7cm}
   \caption{Getting $G$ from $G\,'$: (4/12) canceling an $\{i\}$-dipole}
   \label{fig:CD4}
\end{figure}

\begin{figure}[!h]
   \begin{center}
      \leavevmode
     \includegraphics[width=16cm]{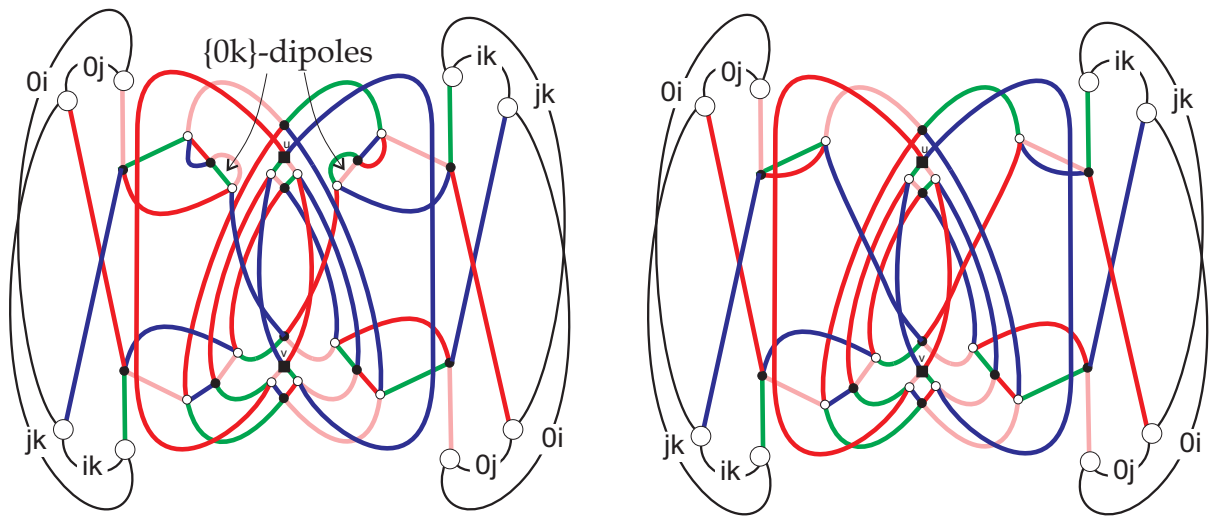}
   \end{center}
   \vspace{-0.7cm}
   \caption{Getting $G$ from $G\,'$: (5/12) canceling two $\{0,k\}$-dipoles}
   \label{fig:CD5}
\end{figure}

\begin{figure}[!h]
   \begin{center}
      \leavevmode
     \includegraphics[width=16cm]{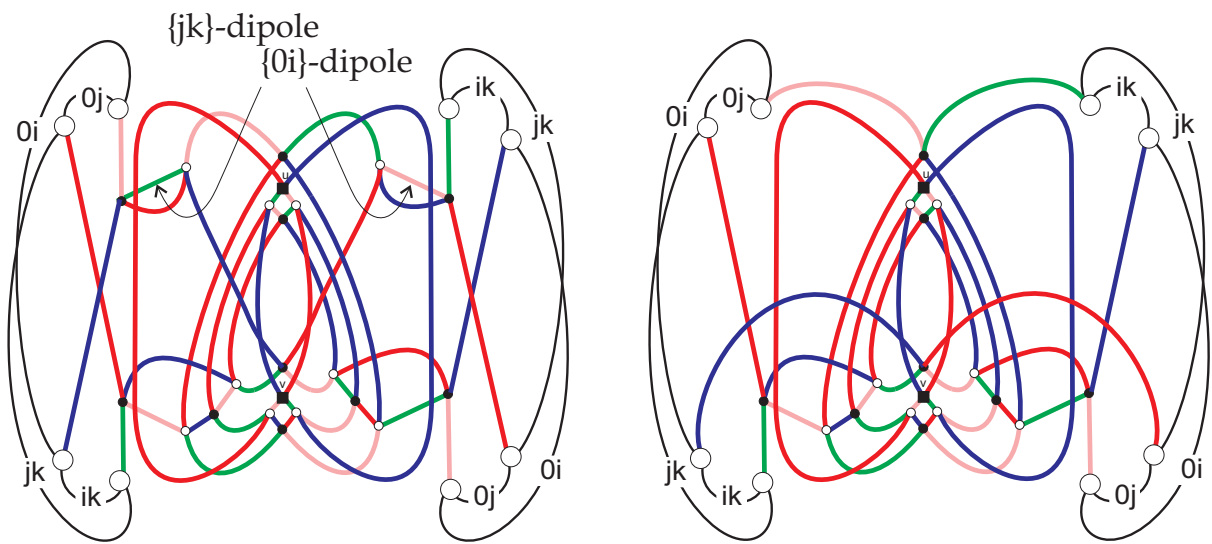}
   \end{center}
   \vspace{-0.7cm}
   \caption{Getting $G$ from $G\,'$: (6/12) canceling $\{0,i\}$-
   and $\{j,k\}$-dipole. No more dipoles}
   \label{fig:CD6}
\end{figure}

\begin{figure}[!h]
   \begin{center}
      \leavevmode
     \includegraphics[width=16cm]{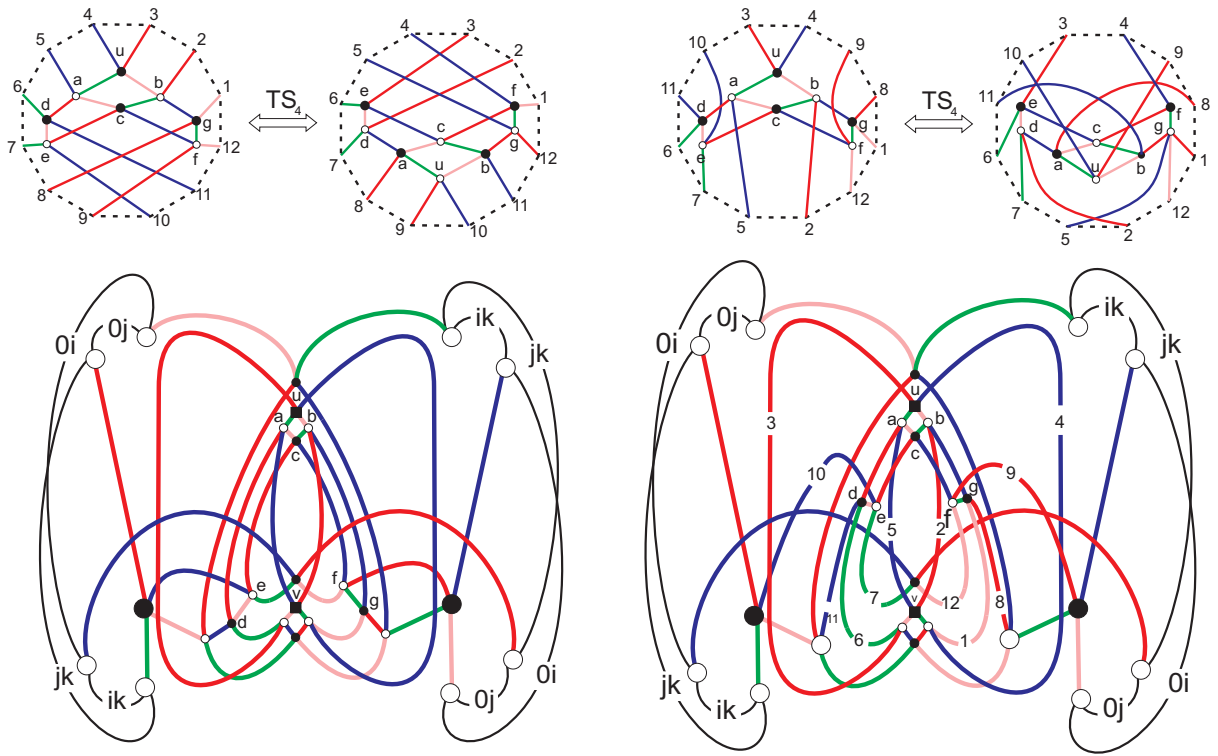}
   \end{center}
   \vspace{-0.7cm}
   \caption{Getting $G$ from $G\,'$: (7/12) aligning move
   $TS_4$ (which factors by dipoles, see pag. 135 of \cite{Lins 1995}}
   \label{fig:CD7}
\end{figure}

\begin{figure}[!h]
   \begin{center}
      \leavevmode
     \includegraphics[width=16cm]{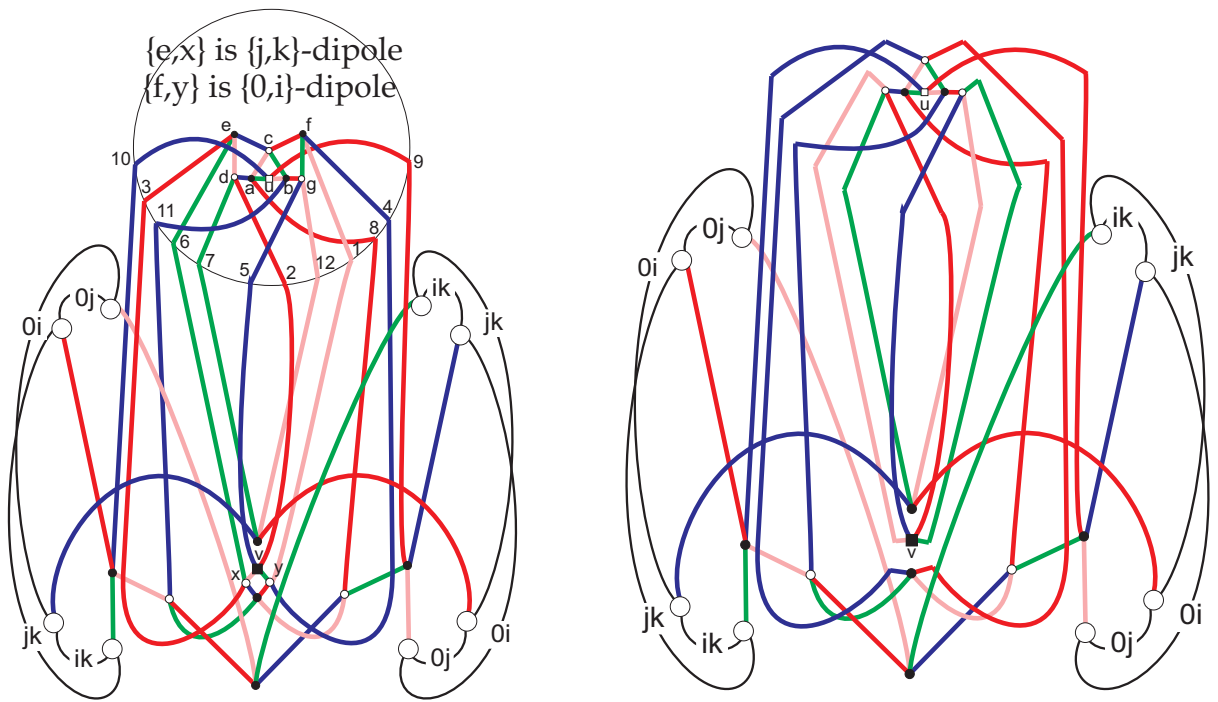}
   \end{center}
   \vspace{-0.7cm}
   \caption{Getting $G$ from $G\,'$: (8/12) canceling $\{0,i\}$-dipole
   and $\{j,k\}$-dipole}
   \label{fig:CD8}
\end{figure}

\begin{figure}[!h]
   \begin{center}
      \leavevmode
     \includegraphics[width=16cm]{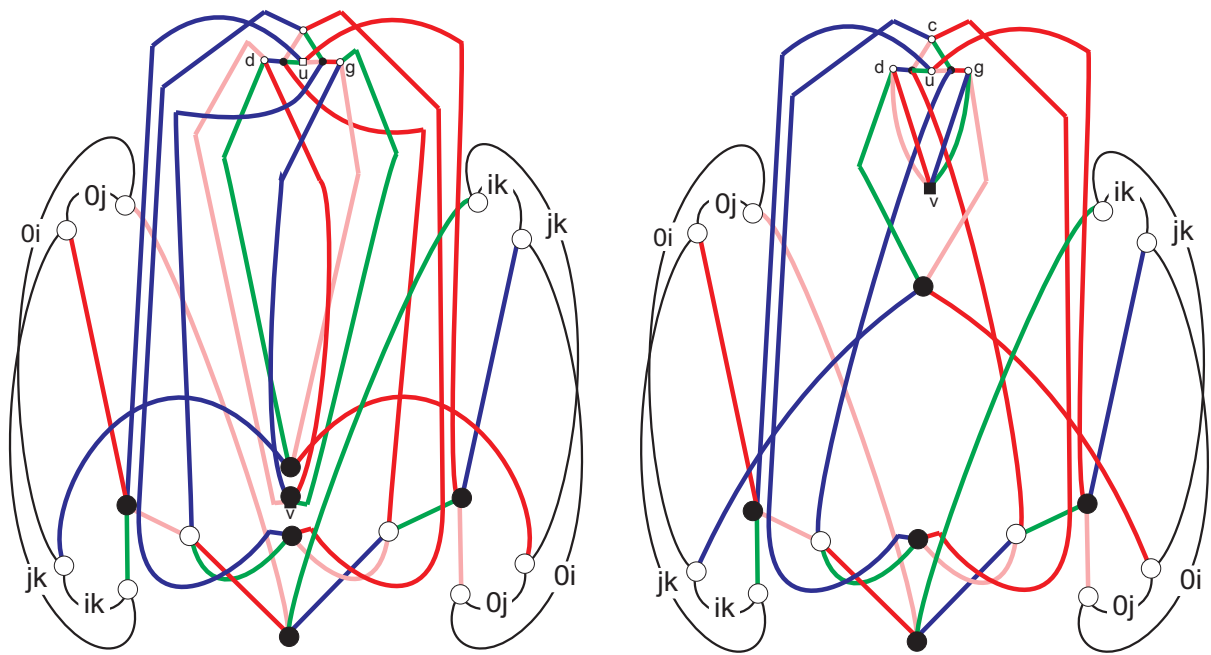}
   \end{center}
   \vspace{-0.7cm}
   \caption{Getting $G$ from $G\,'$: (9/12) aligning the drawing}
   \label{fig:CD9}
\end{figure}

\begin{figure}[!h]
   \begin{center}
      \leavevmode
     \includegraphics[width=16cm]{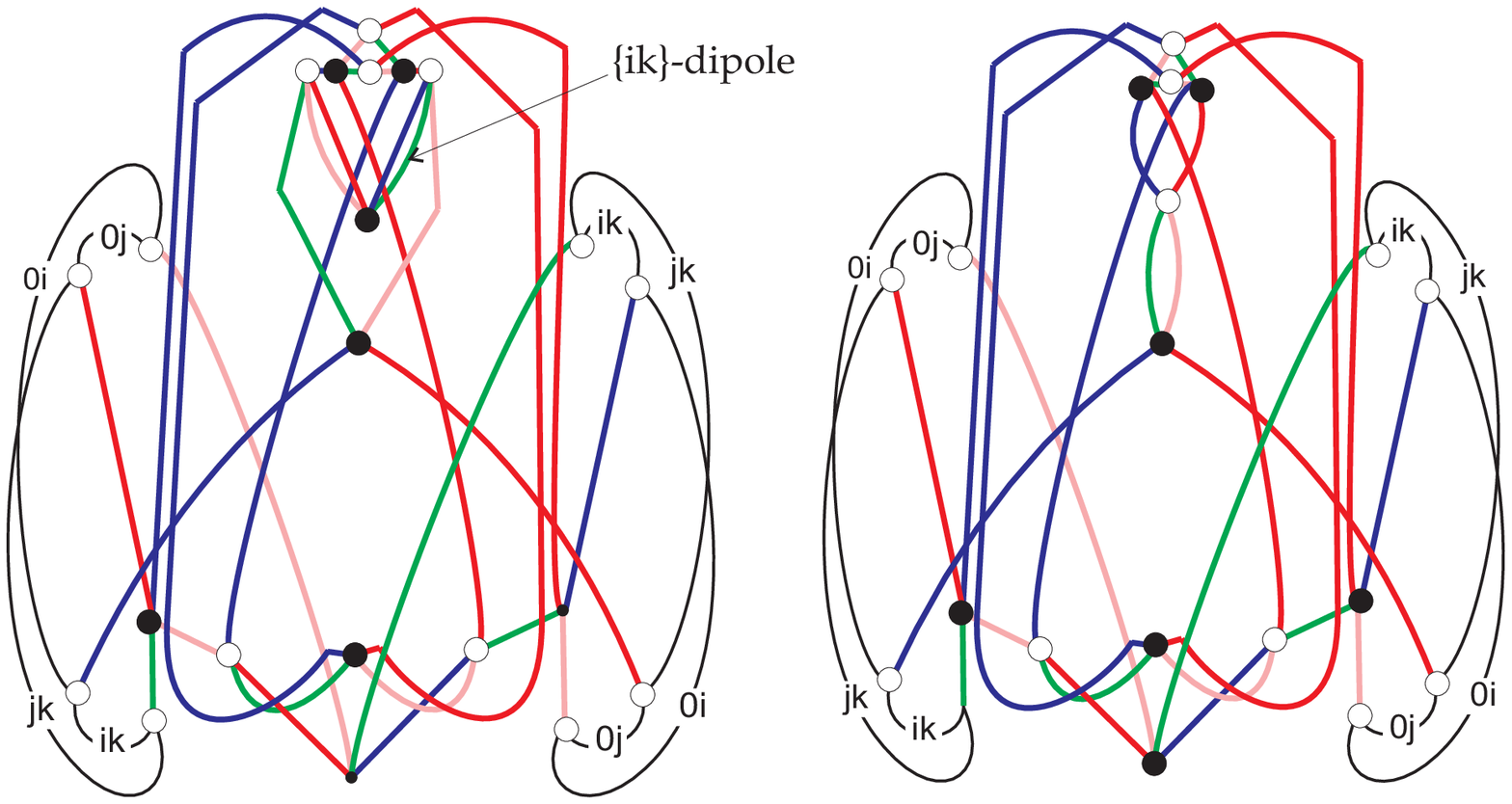}
   \end{center}
   \vspace{-0.7cm}
   \caption{Getting $G$ from $G\,'$: (10/12) canceling $\{0,j\}$-diplole}
   \label{fig:CD10}
\end{figure}

\begin{figure}[!h]
   \begin{center}
      \leavevmode
     \includegraphics[width=16cm]{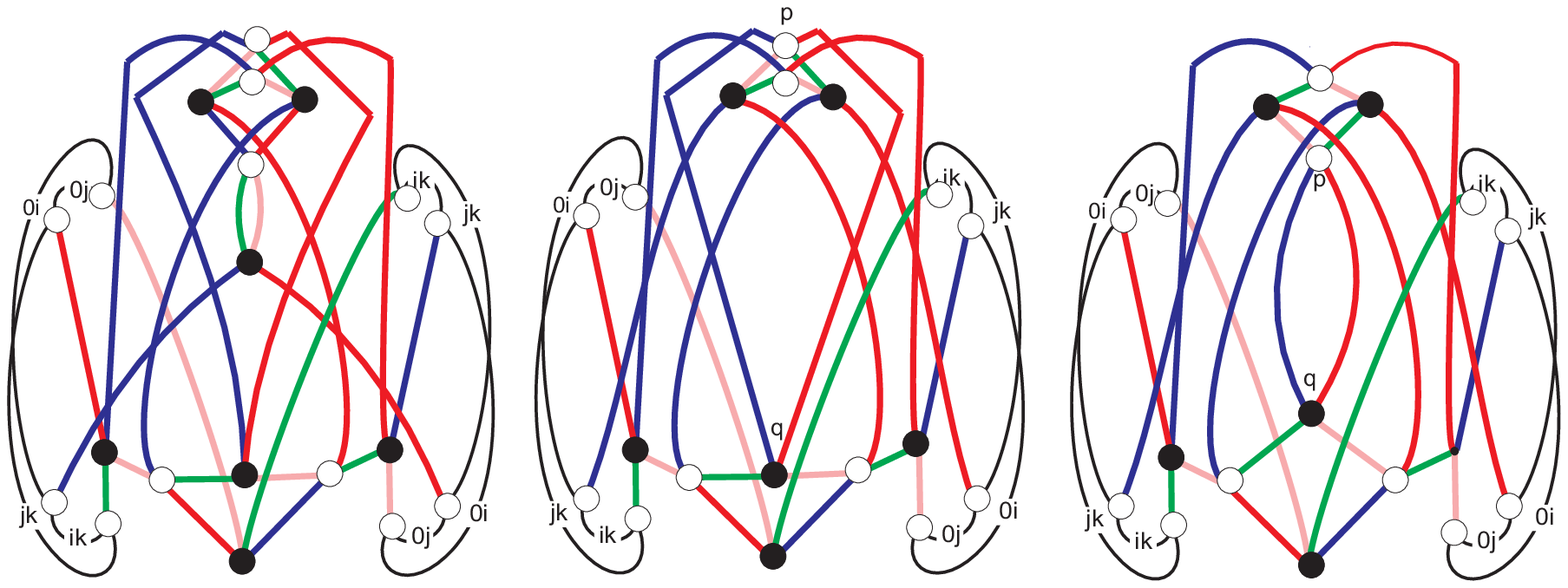}
   \end{center}
   \vspace{-0.7cm}
   \caption{Getting $G$ from $G\,'$: (11/12) canceling $\{0,k\}$-diplole}
   \label{fig:CD11}
\end{figure}

\begin{figure}[!h]
   \begin{center}
      \leavevmode
     \includegraphics[width=16cm]{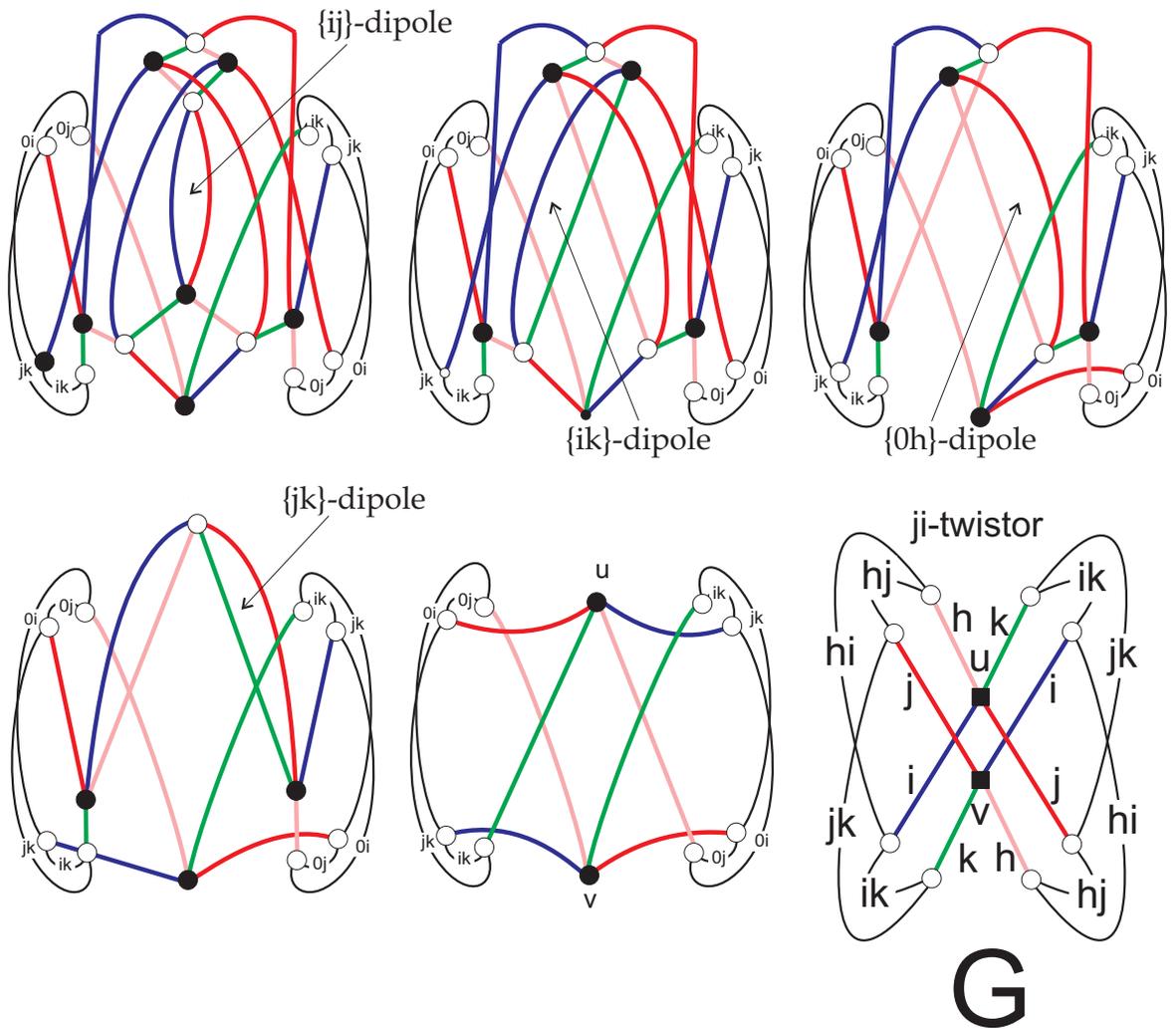}
   \end{center}
   \vspace{-0.7cm}
   \caption{Getting $G$ from $G\,'$: (12/12) $2$-dipole cancelations and
   $u,v$ label interchange}
   \label{fig:CD12}
\end{figure}

\end{document}